\documentstyle[12pt]{article}
\title{ On Local Rank,  Joinings and  Asymptotic Properties of Measure-Preserving Actions}
\author{V.V. Ryzhikov}
\date{}
\def\eps{\varepsilon}
\def\Z{{\bf Z}}
\def\R{{\bf R}}

\def\u{\bigsqcup}
\def\eps{\varepsilon}
\date{vryzh@mail.ru}
\begin{document}

\maketitle
\begin{abstract}    The  note contains a collection of facts and observations 
 around locally rank one actions as well as  constructions 
 connected with some  results by  T. Downarowicz, A.Katok, J.King,  
 F.Parreau, A.A. Prikhodko,  E.Roy, J.Serafin, J-P.Thouvenot et al.  
\end{abstract}
\vspace{3mm}
{\large

\section{ Definitions. Results.  } 
In this note we are speaking
 about connections between rank properties, self-joinings   and certain asymptotic  properties of measure-preserving actions. 
Using the term ``action'' instead of ``transformation'' 
 we hint at the fact that most  of our results
 can be naturally generalized to $\Z^n$- and $\R^n$-actions.  
 Thus, let us recall  necessary definitions for $\Z$-actions only.  

\bf  Local rank. \rm Let $T$ be an automorphism of  a measure space $(X,\mu)$, where $\mu(X)=1$. 
 We denote by  
 $\beta(T)$ its \it  local rank. \rm  It is defined as   maximal  number $\beta$ such that
there exists a sequence of finite partitions of the form 
$$ \xi_j=\{ B_j,  TB_j,    T^2B_j,    \dots,  T^{h_j-1}B_j,C_j^1,\dots, C_j^{m_j}\dots\}$$
for which $\mu(U_j)\to \beta$ and 
any measurable set can be approximated by 
$\xi_j$-measurable ones as ${j\to\infty}$. The sequence of sets  
$$U_j=\bigsqcup_{0\leq k<h_j}T^kB_j$$ is called a~sequence of $\beta$-towers. 
One writes $ \xi_j\to\eps$ and call the sequence
of tower par\-ti\-tions $\{ B_j,  TB_j,    T^2B_j,    \dots,  T^{h_j-1}B_j \}$ {\it approximating}.

\bf Joinings. \rm A {\it self-joining} (of order 2) is defined to be a $T\times T$-invariant
measure $\nu$ on $X\times X$ with the marginals  equal to  $\mu$:
$$\nu(A\times X)=\nu(X\times A)=\mu(A).$$
A joining $\nu$ is called ergodic if the dynamical system
$(T\times T, X\times X, \nu)$ is ergodic.
The measures $\Delta^i= (Id\times T^i)\Delta$ 
defined by the formula
$$
	\Delta^i(A\times B) = \mu(A\cap T^iB) 
$$ 
are referred to as {\it off-diagonals measures\/} (for $i\neq 0$). 
If $T$ is ergodic, then $\Delta^i$
are   ergodic self-joinings. We say that  $T$  has {\it minimal self-joinings}
 of   order 2 (and we write $T\in MSJ(2)$)  if $T$  has no  ergodic
joinings except $\mu^{\otimes 2}=\mu\times\mu$ and $\Delta^i$.

\bf Patial mixing. \rm An automorphism $T$ is called {\it partially mixing}, if
 for for some $\alpha \in (0,1]$ and for all measurable sets $A,B$ 
$$\liminf_i 
\mu(A\cap T^iB)  \geq \alpha \,\mu(A)\mu(B).$$
The maximal value of $\alpha$ satisfying this property for a given $T$ is denoted as $\alpha(T)$, 
and $T$ is called mixing if $\alpha(T)=1$.

\bf { Mildly mixing}.  \rm We say that an automorphism $T$ is {\it mildly mixing} if
 for any set $A$ of positive measure  
$$\limsup_i 
\mu(A\cap T^iA)  < \mu(A) .$$
Actually, $T$ is mildly mixing if it has no rigid factors (except trivial one). 
Let us recall that a \it factor \rm is a restriction of our action onto
some invariant $\sigma$-algebra.  \vspace{2mm}

\bf Partially rigidity. \rm An automorphism  $T$ is said to be {\it partially rigid} if
 for all measurable $A,B$ for some $\rho \in (0,1]$
$$\limsup_i 
\mu(A\cap T^iA)  \geq \rho \mu(A).$$
\vspace{3mm}
The maximal $\rho$ is denoted  by $\rho(T)$ ( $T$ is called rigid as $\rho(T)=1$).

\bf RESULTS. \rm In this note we  present the following results.  

{\bf (1)} Calculations of local rank.  For any $ b\in (0,1)$ there exists a partially mixing action  $\Phi$ of local rank  $\beta (\Phi)=b$; 

\bf Remark.  \rm Let us  recall  A.Katok's  result: \it There exists a weakly mixing $T$ 
such that  the local rank $\beta(T\otimes T)\geq \frac 1 {4}$. \rm
This property of $T$ is proved to be generic~\cite{K}.  
It~can be shown that  $\beta(T\otimes T)\leq \frac 1 4$ and    $Rank(T\otimes T)=\infty$ (see \cite{R1}).
A similar approach to getting a positive local rank  works well for a Cartesian product of two  different actions.
In paper~\cite{KT} J.King and J-P.Thouvenot  state 
as a remark the following assertion: \it  Given $a >0$ there are automorphisms  $S_1,S_2$ such
that $T=S_1\times S_2$ is partially mixing,  and  $\beta(T)>1-a$.\rm

 {\bf (2)} We propose a mildly mixing transformation close to rank one and far from MSJ-class.  
For any $\eps>0$  in the class $\{T:\ \beta(T)>1-\eps \}$ there is a mildly mixing transformation  
with uncountable centralizer and uncountable system of factors. 
 
We provide an infinite  system of factors  for a partially mixing  action  
with the positive local rank via   a modified Katok's construction. 
For any $ b\in (0,\frac{1}{4})$  we find a partially mixing  action $\Phi\otimes\Phi$ of local rank $b$.

{\bf (3)}  A generalization of  King-Thouvenot's theorem to $\Z^n$-actions with ${n>1}$.  
Given a  $\Z^n$-action $\Phi$ satisfying  property $\alpha(\Phi)+\beta(\Phi)>1$,  we show that 
$\Phi$  has   the finite-joining-rank property.\ Its  Markov centralizer  is generated by combinations of action elements with a fixed finite collection of  Markov operators. 

For rank one partially mixing $\Z^2$-actions
 nontrivial self-joinings are  discovered  by T. Downarowicz and J. Serafin \cite{DS}.

{\bf (4)} In contrast to the finite-joining-rank situation and in connection with Proposition 18 of~\cite{PR}, where 
the authors  proved  that a non-rigid  Poisson suspension is not of rank one,  we consider actions of  a positive local rank  assuming the existence of a sequence   of special   joinings tending to $\Delta$.
In this case (close to one discussed in the proof of the mentioned result of~\cite{PR}) 
we prove  the partial rigidity of such  actions.   
This observations, combined with the   proof  by F. Parreau and  E. Roy,  gives: \it 
non-rigid  Poisson suspensions are of zero local rank.  \rm

\bf Remark.  \rm  It's an open question: do there exist actions $T$ and $S$ 
with the same spectral but different rank  properties.
If a Poisson suspension $T_*$ has a positive local rank, then we get a contrast with Gaussian automorphism $G_T$.
  They are of the same spectrum identical to ${{\bf exp}(\sigma_T)}$,  
  but $G_T$ has zero local rank (T. de la Rue).    

{\bf (5)} Banach's problem and Rokhlin's multiple mixing problem are connected within the class of automorphisms
having positive local rank.  So, there is an additional  interest 
to study of multiple mixing properties of locally rank one systems.
In \cite{P} A.A. Prikhodko  introduced new examples of iceberg transformations with a positive local rank.   Jointly with A.A. Prikhodko
we state a {\em positive\/} result:  mixing iceberg transformations are mixing of all orders.

The author thanks Alexander  Prikhodko  for helpful comments.

\section{ A positive local rank for  tensor products of partially mixing systems.  }

\bf Calculating local rank. \rm We explain how to build a pair of 
partially mixing automorphisms $S,\tilde S$ such that for $T=S\times \tilde S$ 
the local rank   $\beta(T)>1-a$. \rm
Let~us consider two rank one constructions for $S$ and $ \tilde S$ with  height sequences $h_j$ and $\tilde h_j=h_j+1$.
Let the corresponding  cut sequences  $r_j$ and $\tilde r_j$ satisfy the condition  $r_j,\tilde r_j \gg (h_j)^3$.  
As in the case of spacer sequences assume that our rank one construction 
have constant (flat) parts and  mixing
stochastic Ornstein's parts.  This   provides $$\beta(T)=(1-\eps)(1-\tilde\eps)>1-a, \ \ \alpha(T)= \eps\tilde\eps.$$ 
Indeed,  we have: 

1. $\alpha(T)= \alpha(S) \alpha(\tilde S)=\eps\tilde\eps.$

2. $\beta(T)\geq (1-\eps)(1-\tilde\eps)$ (we have a sequence of approximating $\beta$-towers of heights $h_j(h_j+1)$ and of measures 
tending to $(1-\eps)(1-\tilde\eps)$).

3. $\beta(T)\leq (1-\eps)(1-\tilde\eps)$.

To prove this we use the following property of a transformation $T$ with $\beta(T)>0$ (see \cite{R1}).    
If $T$ commutes with $R$, then for any $\delta>0$
there is an $m>0$ and a sequence $n(i)$ such that  
$$T^{n(i)}\to (\beta(T)-\delta)R^m +\dots$$
Our notation ``$\dots$'' always means ``something positive''. In the above expansion 
we have  a sum of weighted Markov operators containing some  ``non-$R^m$-part'' 
$$\eps(1-\tilde\eps) \Theta\otimes \tilde P  +  (1-\eps) \tilde\eps P \otimes \Theta + \eps\tilde\eps\Theta\otimes \Theta,$$
where $P,\tilde P$ are some Markov operators and  $\Theta$ is the orthoprojection onto the space of constants  in $L_2(X,\mu)$. 
Setting $R=S\otimes I$,  for arbitrary  $\delta>0$ we get
$$   \beta(T)-\delta  \leq   (1-\eps)(1-\tilde\eps).$$

\bf Modified Katok's examples. \it  For any given $a,b>0$ with $a+2b=1$ there exists an automorphism  $T$ such that 

$\alpha(T)=a^2, \ \beta(T)=b^2$, 	
\\
and $T$ possesses an infinite structure of factors. \rm

For this we combine   Katok's construction (giving positive local rank of $S\otimes S$) with Ornstein's mixing transformation (providing partially mixing).  For instance, consider a rank one construction $S$ with the spacer sequence  
$$ s_j(1), s_j(2),\dots, s_j(q_j), 0, 0, \dots, 0, 1, 1, \dots, 1$$
having $q_j$   zeros and $q_j$  ones with $q_j \gg h_j$ corresponding to the Katok's construction part 
and stochastic $s_j(i)$ (Ornstein's part).
Then for almost all rank one constructions of such kind the generated transformation $S$ 
has the property ${\alpha(S\otimes S)=\beta(S\otimes S)=1/9 }$.    
Now we observe a family of factors $(I\otimes T^n)\cal F$ , where $\cal F$ is  the symmetric factor-algebra  of $S\otimes S$.  
We note also that    for these  constructions one has $\beta(S\odot S)=2\beta(S\otimes S)$.

\bf Mildly mixing, local rank and  an abundance of joinings. \rm
In connection with the question due to J.-P.Thouvenot   
on MSJ property of mildly mixing rank one transformations ($T$ is of Rank 1 iff $\beta(T)=1$)
we remark that 
a little ``rank freedom'' could imply a certain abundance of joinings.

\bf  Examples. \it   Given $\eps >0$ there is a mildly mixing automorphisms  $T$ having
$\beta(T)>1-\eps$, and possessing both an uncountable centralizer and an uncountable structure of factors. \rm

Let us apply King-Thouvenot's idea for infinite products.
We consider partially mixing automorphisms $S_i$ such that 
$$\beta(\bigotimes_{i=1}^{\infty} S_i) =\prod_{i=1}^\infty (1-\alpha(S_i))>1-\eps, $$
$$\alpha(\bigotimes_i^N S_i) =\prod_i^N (\alpha(S_i)). $$
Since all the finite products $\otimes_i^N S_i$ are partially mixing, 
the infinite product has to  be mildly mixing as well. Indeed,
if there exists a rigid function $f$, then all its projections 
to the finite product are rigid, hence, the  projections are  constant functions, so $f$ is constant.

\section{Joinings of partially mixing actions with positive local rank.}
The {\it off-diagonals\/} measures $\Delta^z= (Id\times T^z)\Delta$ 
are defined by the formula
$$\Delta^z(A\times B)=
\mu(A\cap T^zB).$$
If an action  is ergodic, then $\Delta^z$
are   ergodic self-joinings. $T$ 
is called action with {\it minimal self-joinings of   order $2$}  
whenever $T$  has no  ergodic
joinings except $\mu\times\mu$ and $\Delta^z$.

Let us represent a joining $\nu$ in the  form 
 $$\nu(A\times B)=\int_A \nu_x(B) d\mu(x),$$
where $\nu_x$ are conditional measures.
The relatively independent  product $\nu\otimes_X\nu'$ is defined by the equation 
$$\nu\otimes_X\nu'(A\times B)=\int_X \nu_x(A)\nu'_x(B) d\mu(x).$$
A reader familiar with Markov operators in $L_2(\mu)$ can see
that  $\nu\otimes_X\nu'$ corresponds to the operator $P'^\ast P$,
where $(P\chi_A, \chi_B)=\nu(A\times B)$ (here the Markov operator $P$ corresponds to the polymorphism $\nu$). 

\bf  Theorem 3. \it If $\alpha(\Phi)+\beta(\Phi)>1$, then $\Z^n$-action $\Phi$  
has finitely many classes of equivalent ergodic joinings 
($\nu \sim \nu'$ iff $\nu=(Id\times T^z)\nu'$),   
and all ergodic non-trivial self-joinings are the graphs of finite-value maps. \rm

{ \bf Remark. } For  $n=1$ King and Thouvenot proved   minimal self-joinings \cite{KT}.  
Non-trivial joinings could appear for $n>1$
 \cite{DS}. 

\rm
Proof.  For simplicity we consider the case when $\Phi$ is a $\mathbf Z$-action 
generated by a single transformation $T$.  The proof in the general case is similar 
(and we assume that  rank-towers are rectangle). 
Let $\nu\neq \mu^{\otimes 2}$ be an ergodic self-joining.
Denote
$$\eta:=\nu\otimes_X\nu.$$
If 
$$\eta=a\Delta +...,   \qquad  a>0, \eqno (1)$$
then $\nu$  is situated on  a graph of a finite-value map.
We note that the case
$$\eta=a\Delta^z +...,   \qquad  a>0, \quad z\neq 0,$$
is impossible. It means that $\nu$ and $(I\otimes T^z)\nu$ are not disjoint, so $\nu=(I\otimes T^z)\nu$.
This implies $\nu=\mu\otimes\mu$, since $T^z$ is ergodic.

Our aim is to prove (1).   Suppose  $\eta$ and $\Delta$ to be  disjoint, and 
expand $\eta$ in sum 
$$\eta=c_{max} \ \mu^{\otimes 2} +(1-c_{max})\ \tilde\eta,$$
where the number $c_{max}$ is maximal. We see that 
$$\mu^{\otimes 2} \quad  and \quad \tilde\eta \quad  are \ disjoint.\eqno (2)$$
If $c_{max}=1$, then $\eta =\mu\times\mu=\nu$.  So   let  $c_{max}< 1$. 

Given $\eps>0$, $0\leq k\leq \eps h_j$, we define the sets $C^k_j$
 (called columns):
 $$ C^k_j=\u_{i=0}^{h_j-k} T^i T^kE_j\times T^i B_j.
 $$
 For negative $k$ $(-\eps h_j\leq k\leq 0)$, we put
 $$ C^k_j=\u_{i=0}^{h_j+k}T^i B_j\times T^i T^{-k}B_j.
 $$
 Next, we define 
 $$D^\eps_j=\u_{ k=-[\eps h_j]}^{[\eps h_j]} C^k_j,  \qquad  U_j^\eps =\u_{ k=0}^{[\eps h_j]} T^kB_j.,$$
 and for   $\eps > 0$ we have
 $$ \eta(D^\eps_j) >  \eps^2  \beta(T)^2. $$
To see this we use the facts that 
$$U_j^\eps\times U_j^\eps  \ \subset   D^\eps_j, \ \ \ 
 \eta(U_j^\eps\times U_j^\eps) \geq  \mu(U_j^\eps)^ 2,    \ \ \   \mu( U_j^\eps)\approx \beta(T) \eps. $$

Given $\eps>0$   we  obtain 
$$\tilde\eta(D_j^\eps)\geq  d(\eps)>0.$$
For a subsequence $j'$ denoted again by $j$ one has
$$\tilde\eta(\cdot\; |  \ D_j^\eps)\to \eta' \ll \tilde\eta .$$

Inside   the domains $ D_j^\eps$ our joining $\tilde\eta$ is approximated by parts of off-diagonals $\Delta_j^z=\Delta^z( \ |C^z_j)$ (see \cite{R2}): 
$$\lim_j \tilde\eta(\ |D_j^\eps) = \lim_j \sum_{|z|<\eps h_j} a^z_j\Delta^z_j,$$
where $a^z_j=\tilde\eta(C^z_j|D_j^\eps)$. 
Note that  $\sum_{|z|<M} a^z_j$ tends to zero  for fixed $M$ and $j\to\infty$ (we suppose now that $\tilde \eta\perp \Delta^z$).
If $z$ is large,  then
$$\Delta^{z}\approx  \alpha(\Phi)\mu^{\otimes 2} +...\ . $$
$$\Delta^{z}_j\approx  [\alpha(\Phi) -\eps\beta(\Phi) - (1-\beta(\Phi) )]\mu^{\otimes 2} +...\  $$
whenever $|z|<\eps h_j$.
So  $\tilde\eta$ contains  a part $$\lim_j \sum_z a^z_j\Delta^z_j=c\mu^{\otimes 2}+\dots, $$
where $c\geq \alpha(\Phi) +\beta(\Phi) - 1-\eps \beta(\Phi) >0$.

Hence, $\mu^{\otimes 2}$ and $\tilde\eta$ are not disjoint, but this  contradicts (2)! So our assumption on 
the disjointness  of $\eta$ and $\Delta$ is  false. 
We have proved 
$\eta \gg \Delta$ (in fact $\eta =\frac{1}{N} \Delta+\dots$), so  $\nu$ is a graph of a finite-value map.

A composition of two  finite-value maps is   a finite-value map. This implies that a relative product
$\nu\otimes_X\nu'$  of two non-trivial joinings  is supported on a graph,  hence,  $\nu\otimes_X\nu'$  
and $\mu^{\otimes 2}$
are disjoint.

If we have a collection $\nu_1, \dots, n_k\neq \mu^{\otimes 2}$ of ergodic self-joinings, and $k> [(\alpha +\beta - 1)^{-1}]+1$,
then there exist different $m,n$ ($1\leq m,n \leq k$) such that  
$$\limsup_j \nu_m\otimes_X \nu_n (U_j^\eps\times U_j^\eps) > 0.$$
The mentioned  approximation shows that if $\nu_m\otimes_X \nu_n$ has no $ \Delta^z $ as a component, 
then it must contain as a component $\mu^{\otimes 2}$, but this is impossible.
Thus,
$$\nu_m\otimes_X \nu_n =a\Delta^z +\dots\ .$$
This means that  $\nu_m$ and $\nu_n$ are equivalent.

\section{Local Rank and Partial Rigidity}
{\bf Remark.  }\it  Consider a locally rank one ergodic $Z^n$-action $T$. 
If there is a  joining $\nu\neq \Delta^z$ ($P\neq I$) for $T$ 
with $\nu\otimes_X\nu \gg \Delta$ (respectively, $P^*P=cI+\dots$), 
then $T$ is partially rigid. \rm

Indeed,  if $T^{z(i)}\to a P +\dots$, $P\neq T^z$,  then there is  $w(j)\to\infty$ such that 
$$T^{w(j)}\to a^2 P^*P +\dots = c'I+\dots, \ c'>0.$$   So our action is partially rigid.

In a recent preprint \cite{PR}  F. Parreau and E. Roy  
proposed an interesting approach to the study of rigidity and  rank of  Poisson suspensions
 using  a certain natural sequence of  joinings $\nu_i\to \Delta$.
Their work stimulated the writing of the following theorem.

\bf  Theorem 4. \it Let $T$ be a locally rank one ergodic action. 
If it has a sequence of self-joinings $\nu_i$ with 
$\nu_i\otimes_X\nu_i\perp \Delta$ and  $\nu_i\to \Delta$ then $T$ is partially rigid. \rm

Let us consider a self-joining $\nu$ with a joining $\eta$ corresponding to the operator $P^*P$, 
where Markov operator  $P$ corresponds to  $\nu$.  Denote 
$\eta=\nu\otimes_X\nu$ and consider ${U^\eps_j}$ and  $D_j^\eps$ (as well as $\Delta^z_j$) 
introduced in the proof of Theorem 3.
Since for $\eps>0$
$$ \liminf_j\eta( U^\eps_j\times U^\eps_j)= \liminf_j \left(P^*P\chi_{U^\eps_j}\\ ,\  \chi_{U^\eps_j}\right) \geq c\geq (\eps\beta(T))^2,$$ we get
$$\eta(D_j^\eps)\to d>0, \ \
\eta(\ | D_j^\eps)\to \eta',     \ \  \eta=d\eta' +... . $$
Applying again approximation arguments  we get 
$$\sum_z a_j^z \Delta^z_j  \to \eta' .$$

Now let us turn to our sequence $\nu_i\to \Delta$.  We have $P_i\to_w  I$, where $P_i$ corresponds to the joining $\nu_i$. Since  $I$ is unitary we see that
$P_i\to_s  I$ and $P^*_iP_i\to_s I$ (we~use the symbol $\to_s$ for strong convergence).  
Thus, $\eta_i\to \Delta$, and 
the corresponding sequence $\eta'_i\to \Delta$ as well ( $ \Delta$ is ergodic self-joining!).
We get
 $$\sum_z a_{j(i)}^z \Delta^z_j  \to \Delta. $$
 Choice lemma \cite{R2}  and the possibility to have  $a_{j(i)}^0\to 0$ (as a consequence of $\nu_i\otimes_X\nu_i\perp \Delta$) 
 together provide the  partial rigidity of the action: in fact, there is a sequence ${z_i\to \infty}$ such that
$$ \Delta^{z_i}_{j(i)}  \to \Delta ,  \  \Delta^{z_i}  \to \beta(T)(1-\eps)\Delta +... , \ T^{z_i}  \to \beta(T)(1-\eps) I +\dots . $$
Our action $T$ is partially rigid.

\section{ Multiple Mixing of Iceberg  Transformations}
We know that the condition $\beta(\Phi)>2^{-n}$ for  mixing $\Z^n$-action $\Phi$ implies  multiple mixing.   
If a  mixing $T$ of positive local rank
hasn't multiple mixing, then due to B.\,Host et al.\ 
we can say that an automorphism is found with absolutely continuous   spectrum
of finite multiplicity.   Thus, within  the  class of  mixing transformations 
of a positive local rank a connection is observed between the spectral  problem on the existence of transformations with absolutely continuous spectra of finite multiplicity and  Rokhlin's multiple mixing problem.   
In \cite{P} new examples of transformations with positive local rank are introduced.   
Although, instead of new spectral phenomena we got the following result of positive kind. 
 
\bf Theorem 5\rm  (jointly with A.A.Prikhodko). \it  Mixing iceberg  transformations are mixing of all orders. \rm

In fact, this theorem follows from certain  properties of the construction \cite{P}  (there is some analogy with  finite rank actions), 
 and the theorem on multiple  mixing of systems with $D$-property \cite{R}.

}

 \end{document}